\documentclass[12pt,a4paper]{article}
\usepackage{amsmath,amssymb,latexsym,amsthm}
\usepackage[english]{babel}

\date{}

\title{On geodesic mappings of  manifolds with affine connection}
\author{Josef Mike\v s
and Irena Hinterleitner
\thanks{Supported by the Council of the Czech Government MSM 6198959214.}}

\begin{document}
\def\be#1{\begin{equation}\label{#1} }
\newcommand{\ee}{\end{equation}}
\newcommand{\bea}{\begin{eqnarray}}
\newcommand{\eea}{\end{eqnarray}}
\newcommand{\bc}{\begin{center}}
\newcommand{\ec}{\end{center}}
\def\journal#1#2#3#4{{#1} {\bf #2}, #3 (#4)}
\newtheorem{proposition}{Proposition}
\newtheorem{lemma}{Lemma}
\newtheorem{theorem}{Theorem}
\newtheorem{corollary}{Corollary}
\newtheorem{Remark}{Remark}
\def\eps{\varepsilon}
\def\GG{\Gamma}
\def\g{\gamma}
\def\phi{\varphi}
\def\epsilon{\varepsilon}
\def\ep{\varepsilon}
\def\s{\sigma}
\def\t{\tau}
\def\r{\varrho}
\def\e{{\rm e}}
\def\ds{\displaystyle }
\def\fr#1#2{{\ds\frac{#1}{\ds #2}}}
\def\noi{\noindent}
\def\cotan{{\,\rm cotan \,}}
\def\const{{\,\rm const \,}}
\def\sign{{\,\rm sign \,}}
\def\cc#1={\cite{#1}}
\def\R{\hbox{$\mathbb{R}$}}
\def\RR{\hbox{$\mathbb{R}^2$}}
\def\Rn{$\mathbb{R}^n$}
\def\L{\hbox{$\cal{L}$}}
\def\a{\alpha}
\def\b{\beta}
\def\G{\Gamma}
\def\d{\delta}
\def\s{\sigma}
\def\nad#1#2{\buildrel{#1} \over{#2}\!\!\strut}
\def\pod#1#2{\mathrel{\mathop{#2}\limits_{#1}}\strut}
\def\cc#1={\cite{#1}}
\def\ttt{\hbox{\,.\,.\,.\,}}
\def\an{\hbox{$A_n$}}
\def\ann{\hbox{$\bar A_n$}}
\def\vnn{\hbox{$\bar V_n$}}
\def\annn{\hbox{$\tilde A_n$}}
\maketitle

\begin{abstract}
In this paper we prove that all  manifolds with affine connection are globally projectively
equivalent to some space with equiaffine connection (equiaffine manifold). These manifolds are characterised by a symmetric Ricci tensor.
\end{abstract}

\section{Introduction}%

Diffeomorphisms between two spaces \an\ and \ann\ with affine connection
are called {\it geodesic}, if any geodesic of \an\ is mapped to a
geodesic of \ann, see for example \cite{Eastwood-Matveev} -- \cite{Thomas},
etc.

{\it Affine} (or {\it trivial geodesic\/}) {\it mappings} are
mappings which preserve canonical parameters of geodesics. Many
papers are devoted to the metrizability or projective metrizability
of spaces with affine connection. Under metrizability of a space
\an\ we understand the existence of a metric $g$, which generates
the affine connection $\nabla$, such that $\nabla$ is the
Levi-Civita connection of $g$ (for which $\nabla g=0$), see
\cite{Eastwood-Matveev,Eisenhart1926,Eisenhart2000,Mikes-Kiosak-Vanzurova}.
On the other hand, the problem of metrizability, respectively
projective metri\-zability, is equivalent to that of affine,
respectively geodesic, mappings of a space \an\ with affine
connection onto (pseudo-) Riemannian spaces.

{\it Equiaffine spaces}, characterised by the symmetry of the Ricci tensor, are spaces in which the volume of an $n$ dimensional parallelepiped is invariant under parallel transport, play an important role in the theory of geodesic
mappings. J.~Mike\v s and V.E.~Berezovski \cite{Mikes-Berezovski,Mikes-Kiosak-Vanzurova}
found fundamental equations of geodesic mappings of equiaffine
spaces onto (pseudo-) Riemannian spaces in the form of systems of
linear partial differential equations of Cauchy type in terms of covariant derivatives.
These results were used for further studies by M.G.~Eastwood and V.~Mat\-ve\-ev \cite{Eastwood-Matveev}.

In the paper \cite{Mikes-Hinterleitner-Kiosak-08} by I.~Hinterleitner, J.~Mike\v s and V.A.~Kiosak it was proved that any
space with affine connection is locally projective equiaffine. We
show that these properties hold globally, i.e.~any arbitrary manifold with affine
connection globally admits geodesic mappings onto some equiaffine
manifold.

For this reason the solution of the problem of the projective metrizability of a
space~\an\ (or equivalently of geodesic mappings of \an\ onto
(pseudo-) Riemannian manifolds \vnn)  can be realized as geodesic
mapping of the equiaffine space \annn, which is projectively
equivalent to the given space \an.

\section{Main properties of geodesic mappings}

Let \an\ = $(M,\nabla)$ be a manifold $M$ with affine connection
$\nabla$. We suppose that $\nabla$ is torsion-free, i.e.
$\nabla_XY-\nabla_YX=[X,Y]$, where $[\:,\:]$ is the Lie bracket:
$[X,Y]f=X(Yf)-Y(Xf)$ for $f\in {\cal F}(M)$.
Here and after $X,Y, \ttt$ are tangent vectors.

The curvature $R$ of a manifold \an\
 is a tensor field of type $(1,3)$ defined by
\be{f0}
R(X,Y)Z=\nabla_X(\nabla_YZ)-\nabla_Y(\nabla_XZ)-\nabla_{[X,Y]}Z,
\ee
called sometimes also the Riemannian tensor of the connection.

We can introduce the Ricci tensor $Ric$ of type $(0,2)$ as a trace of a linear map, namely:
$$Ric(X,Y)={\rm trace}\{V\mapsto R(X,V)Y\}.$$

A manifold \an\ with a torsion-free (symmetric) affine connection is
called {\it equiaffine manifold}, if the Ricci tensor is symmetric , i.e. (\cite{Mikes-Kiosak-Vanzurova,Norden,Sinyukov})
$$Ric(X,Y)=Ric(Y,X).$$

It is known \cite{Mikes-Kiosak-Vanzurova,Norden,Sinyukov} that (pseudo-) Riemannian spaces are equiaffine manifolds. 

\bigskip
A manifold $A_n=(M,\nabla)$ with affine connection  admits a geodesic mapping
onto $\bar A_n=(M,\bar\nabla)$, if and only if the {\it Levi-Civita equation} holds
\cite{Eastwood-Matveev, Eisenhart2000, Mikes-Kiosak-Vanzurova, Sinyukov}:
\begin{equation}\label{f1}
\bar\nabla_XY = \nabla_XY + \psi(X)Y+ \psi(Y)X,
\end{equation}
where $\nabla$ and $\bar\nabla$ are the affine connections of \an\
and \ann\, respectively,
and $\psi$~is a linear form.

Geodesic mappings with $\psi\equiv 0$ are {\it trivial} or
{\it affine}.

\medskip
After applying formula (\ref{f0}) for the
curvature tensor and expression (\ref{f1}) to a geodesic mapping
\an\ $\rightarrow$ \ann\ we found a relationship between the
curvature tensors $R$ and $\bar R$ of \an\ and \ann:
\be{f2}
\begin{array}{r}
\bar R(X,Y)Z =R(X,Y)Z\hfill \\[2mm]
\qquad\qquad\quad +(\psi(X,Y)-\psi(Y,X))\,Z+
\psi(X,Z)\,Y- \psi(Y,Z)\,X,
\end{array}
\ee
where
\be{f3}
\psi(X,Y)=\nabla_Y\psi(X)-\psi(X)\psi(Y).
\end{equation}
By contraction of (\ref{f2}) we obtain the following relation for
the Ricci tensors
\begin{equation}\label{f4}
\bar Ric(X,Y)=Ric(X,Y)+n\,\psi(X,Y)-\psi(Y,X).
\end{equation}

\section{Main results}

Suppose \an\ $\in C^1$,
i.e. the components of the affine connection $\nabla$ of \an\  are  functions of type $C^1$ on all charts of the manifold $M$.

In our paper \cite{Mikes-Hinterleitner-Kiosak-08} we proved that any manifold \an\ with affine connection
is locally projectively equiaffine.

The following theorem holds generally:\\
{\bf Theorem}
{\it All manifolds \an\ with affine connection are projectively equi\-affine.
}\\[2mm]
\noi{\bf Remark.}
In other words, an arbitrary manifold \an\ $\in C^1$ admits a global geodesic mapping onto an equiaffine manifold \ann, and moreower \ann\ $\in C^1$.
\begin{proof} We prove the  existence of a geodesic mapping
of \an\ = $(M, \nabla)\in C^1$ onto an equiaffine manifold \ann\ =
$(M,\bar\nabla)$.

It is known that on the manifold $M$ exists globally a metric, say $\tilde g$. Now we can construct the metric $\tilde g$ so that
$\tilde g\in C^2$, i.e. the components $\tilde g_{ij}$ of $\tilde g$  in a coordinate domain of $M$ are functions of type $C^2$.
We denote by $\tilde\nabla$ the Levi-Civita connection of $\tilde g$.

 We  construct the one-form $\psi$ in the following way
\be{f5}
\psi(X)=-\frac{1}{n+1}\,{\rm trace}(Y\mapsto \nabla_XY-\tilde\nabla_XY).
\ee
With the help of formula \eqref{f1} applied to $\psi$
we construct globally the affine connection $\bar \nabla$ on $M$. It is evident that
\an\ is geodesically mapped onto \ann\ =
$(M,\bar \nabla)$, and, evidently, \ann\ $\in C^1$.

Now we prove that \ann\ is equiaffine. It is sufficient to calculate
that
$$
\bar Ric(X,Y)=\bar Ric(Y,X),
$$
where $\bar Ric$ is the Ricci tensor of \ann.

From (\ref{f4}) follows
$$
Ric(X,Y)-Ric(Y,X)+(n+1)(\psi(X,Y)-\psi(Y,X))=0,
$$
where $Ric$ is the Ricci tensor of \an. Further we use relation
(\ref{f3}):
$$
Ric(X,Y)-Ric(Y,X)+(n+1)(\nabla_Y\psi(X)-\nabla_X\psi(Y))=0.
$$
After using (\ref{f3}) and (\ref{f5}) we can see that this equation
holds identically.
\end{proof}

\noi{\bf Remark.}
The equiaffine connection $\bar \nabla$ constructed in this way is constructed explicitly from the original connection $\nabla$.

\medskip
Similarly it could be shown that a manifold with projective connection admits a geodesic mapping onto a manifold with equiaffine connection.


\end{document}